\RequirePackage{fix-cm}
\documentclass[pdflatex,sn-mathphys-num]{sn-jnl}

\usepackage{amsmath,amssymb,amsfonts,amsthm,mathtools}
\usepackage{enumitem}
\usepackage{xcolor}

\makeatletter
\@ifundefined{th@thmstyleone}{%
\newtheoremstyle{thmstyleone}{18pt plus2pt minus1pt}{18pt plus2pt minus1pt}%
  {\small\itshape}{0pt}{\small\bfseries}{}{.5em}%
  {\thmname{#1}\thmnumber{\@ifnotempty{#1}{ }\@upn{#2}}%
   \thmnote{ {\the\thm@notefont(#3)}}}%
\newtheoremstyle{thmstyletwo}{18pt plus2pt minus1pt}{18pt plus2pt minus1pt}%
  {\small\normalfont}{0pt}{\small\itshape}{}{.5em}%
  {\thmname{#1}\thmnumber{\@ifnotempty{#1}{ }{#2}}%
   \thmnote{ {\the\thm@notefont(#3)}}}%
\newtheoremstyle{thmstylethree}{18pt plus2pt minus1pt}{18pt plus2pt minus1pt}%
  {\small\normalfont}{0pt}{\small\bfseries}{}{.5em}%
  {\thmname{#1}\thmnumber{\@ifnotempty{#1}{ }\@upn{#2}}%
   \thmnote{ {\the\thm@notefont(#3)}}}%
}{}
\makeatother

\theoremstyle{thmstyleone}
\newtheorem{theorem}{Theorem}[section]
\newtheorem{proposition}[theorem]{Proposition}
\newtheorem{lemma}[theorem]{Lemma}
\newtheorem{corollary}[theorem]{Corollary}
\theoremstyle{thmstylethree}
\newtheorem{definition}[theorem]{Definition}

\newtheorem{remark}[theorem]{Remark}
\newtheorem{question}[theorem]{Question}

\newcommand{\Topz}{\mathbf{Top}_0}
\newcommand{\Topo}{\mathbf{Top}_1}
\newcommand{\Sob}{\mathbf{Sob}}
\newcommand{\HH}{H}
\newcommand{\DetH}{\mathbf{HD}}
\newcommand{\cl}{\operatorname{cl}}
\newcommand{\Irr}{\operatorname{Irr}}
\newcommand{\OO}{\mathcal O}
\newcommand{\SSS}{\mathcal S}
\newcommand{\FFF}{\mathcal F}
\newcommand{\CCC}{\mathcal C}
\newcommand{\DDD}{\mathcal D}
\newcommand{\RRR}{\mathcal R}
\newcommand{\PPP}{\mathcal P}
\newcommand{\RDset}{\mathcal{RD}}
\newcommand{\WDset}{\mathcal{WD}}
\newcommand{\hdop}[1]{d_{#1}}
\newcommand{\hdH}{d_H}
\newcommand{\HdH}{\operatorname{Hd}_{H}}
\newcommand{\ptop}{\mathcal E}
\newcommand{\doi}[1]{\url{https://doi.org/#1}}

\setlist[enumerate]{label=\textup{(\arabic*)},leftmargin=*,itemsep=2pt,topsep=4pt}
\hypersetup{hidelinks,bookmarksdepth=2,
  pdftitle={Subset-system derived topologies and coreflective categories of T0-spaces},
  pdfauthor={Xiaoquan Xu, Xinjie Huang, Yi Yang},
  pdfkeywords={subset system, derived topology, coreflective category, Cartesian closed category, categorical topology, domain theory},
  pdfsubject={Applied Categorical Structures submission}}
\begin{document}

\title[Subset-system derived topologies]{Subset-system derived topologies and coreflective categories of $T_0$-spaces}

\author*[1]{\fnm{Xiaoquan} \sur{Xu}}\email{xiqxu2002@163.com}
\author[1]{\fnm{Xinjie} \sur{Huang}}\email{2418782429@qq.com}
\author[2]{\fnm{Yi} \sur{Yang}}\email{1354848335@qq.com}

\affil*[1]{\orgdiv{School of Mathematics and Statistics},
  \orgname{Minnan Normal University},
  \orgaddress{\city{Zhangzhou}, \postcode{363000}, \country{China}}}
\affil[2]{\orgdiv{School of Mathematics and Statistics},
  \orgname{Jiangxi Normal University},
  \orgaddress{\city{Nanchang}, \postcode{330000}, \country{China}}}

\abstract{Given a subset system $H$ on $T_0$-spaces, we associate with each topology a closure-testing derived topology and study the fixed points of this operator.  The $\HH$-determined topologies above a given topology are closed under intersections and hence determine a least fixed-point hull.  We prove that this hull is obtained by transfinite iteration, is functorial, and defines a right adjoint to the inclusion of $\HH$-determined spaces into $\mathbf{Top}_0$.  Consequently, the corresponding full subcategory is coreflective and complete.  We isolate restriction-stability and open-locality hypotheses under which determinedness is inherited by closed and open subspaces, and we establish a general non-reflectivity criterion.  If $H$ is contained in the directed-subset system, the category of $\HH$-determined spaces is Cartesian closed; its exponentials are obtained by coreflecting the pointwise convergence topology on continuous maps.  The construction subsumes familiar monotone- and irreducible-determined cases while separating the formal consequences of the subset-system axioms from the additional hypotheses needed for heredity and exponentials.  It therefore provides a uniform categorical mechanism for producing complete Cartesian closed categories of $T_0$-spaces with explicit higher-order function spaces, as required in denotational semantics.}

\keywords{subset system, derived topology, coreflective category,
  Cartesian closed category, categorical topology, domain theory}
\pacs[2020 Mathematics Subject Classification]{54B30, 18F60, 06B35, 54D10}

\maketitle

\section{Introduction}

Domain theory provides order-theoretic models of computation in which directed suprema describe approximation and the Scott topology records observable information; see, for example, \cite{Abramsky-Jung-1994,Gierz-et-al-2003,Goubault-2013}. From the viewpoint of denotational semantics, a category of semantic spaces is particularly useful when it supports products and higher-order function spaces; Cartesian closed categories provide exactly this structure through exponentials. A parallel line of work develops the order-theoretic core of domain theory directly on $T_0$-spaces. Monotone determined spaces are a basic example: their topology is detected by directed subsets, and they form a Cartesian closed category extending the Scott-topological setting; see \cite{Erne-2009,Luo-Xu-2017,Yu-Kou-2015}. The question pursued here is whether this categorical infrastructure can be generated uniformly from a functorially chosen class of distinguished subsets, rather than constructed separately for each convergence or approximation notion.

Many constructions in this area describe a topology by specifying which distinguished subsets are allowed to detect open sets.  Replacing directed subsets by irreducible subsets leads to irreducible determined spaces.  In \cite{Xu-Yang-Luo-2026}, the corresponding derived topology, its transfinite idempotent hull, hereditary properties and categorical structure were studied for the irreducible subset system.  The purpose of the present paper is not to replace the irreducible system there by a formal parameter.  Instead, we isolate the fixed-point mechanism behind these constructions and determine what follows for an arbitrary subset system
\[
   \HH:\Topz\longrightarrow\mathbf{Set}.
\]
The notion of subset system, inspired by the uniform closure constructions of Wright, Wagner and Thatcher \cite{Wright-Wagner-Thatcher-1978}, was used in \cite{Xu-2021} to organize several sobriety-type notions.  Chen and Xu \cite{Chen-Xu-2019} also studied topological spaces determined by subset systems through an induced-operator framework.  Here the subset system is used in a different way: from a $T_0$-topology $\tau$ we form a closure-testing operator $\hdH(\tau)$.  This operator is monotone and inflationary on the poset of $T_0$-topologies on the fixed underlying set, but its idempotence is not part of the subset-system axioms.  The central object is therefore its least fixed point above $\tau$.

A second recent uniform approach is due to Chen, Kou and Lyu \cite{Chen-Kou-Lyu-2026}.  Their construction starts from convergence classes and $\mathbb C$-generated spaces; in particular, it recovers monotone determined spaces as a Cartesian closed coreflective category and treats directed, chain and monotone sequential convergence classes in a common framework.  The present construction organizes the data differently.  It begins with a closure-testing operator on the lattice of $T_0$-topologies: the distinguished family $\HH(X)$ depends functorially on the current space, candidate open sets are tested against closures of $\HH$-sets, and the one-step operator may require transfinite iteration before reaching a fixed point.  Thus the two approaches overlap in the monotone-determined case, but the fixed-point completion, its ordinal length and the associated heredity questions are intrinsic to the closure-testing operator studied here.

Our main results form three categorical levels.  First, the $\HH$-determined topologies above a given $T_0$-topology are closed under intersections.  Hence there is a least fixed-point hull $\HdH(\tau)$, obtained by transfinite iteration of $\hdH$; the assignment $\tau\mapsto\HdH(\tau)$ is a functorial closure operator on the poset of $T_0$-topologies on a fixed set.  Second, this hull defines a coreflector: the full category $\DetH$ of $\HH$-determined spaces is coreflective in $\Topz$ and therefore complete.  We also obtain retract stability, a general non-reflectivity criterion, and restriction and localization hypotheses that ensure closed- and open-heredity.  Third, if $\HH\leq\DDD$, then $\DetH$ is Cartesian closed.  The exponential is obtained by applying the coreflector to the pointwise convergence topology on the set of continuous maps.  Thus the same fixed-point mechanism produces the completion, the adjunction and explicit higher-order function spaces, which is the categorical structure most directly relevant to denotational semantics.

The paper is organized around these structural points.  Section~2 introduces subset-system derived topologies and records the basic order-theoretic relations among them.  Section~3 develops the fixed-point completion, its transfinite description and functoriality.  Section~4 turns the completion into a coreflector and derives completeness, retraction stability and non-reflectivity criteria.  Section~5 treats subspaces under explicit stability assumptions on $\HH$.  Section~6 constructs exponentials and establishes Cartesian closedness when $\HH\leq\DDD$.  The final section records questions suggested by the general framework.

\section{Subset-system derived topologies}

We recall only the notation needed below.  For a $T_0$-space $X$, the specialization order is denoted by $\leq_X$; thus $x\leq_X y$ iff $x\in\cl_X\{y\}$.  We write $\alpha(X)$ for the Alexandroff topology of upper sets of this order.  A nonempty subset $D$ is \emph{directed} if every finite subset of $D$ has an upper bound in $D$.  The family of directed subsets of $X$ is denoted by $\DDD(X)$, and $\Irr(X)$ denotes the family of nonempty irreducible subsets of $X$.

\begin{definition}\label{def:subset-system}
A covariant functor $\HH:\Topz\to\mathbf{Set}$ is a \emph{subset system} if, for every $T_0$-space $X$,
\[
   \SSS(X)\subseteq \HH(X)\subseteq 2^X,
\]
where $\SSS(X)=\{\{x\}:x\in X\}$, and for every continuous map $f:X\to Y$ one has
\[
   \HH(f)(A)=f(A)\in\HH(Y)~(A\in\HH(X)).
\]
If $\HH_1(X)\subseteq\HH_2(X)$ for every $X$, we write $\HH_1\leq\HH_2$.
\end{definition}

Standard examples are the systems $\SSS,\FFF,\CCC,\DDD,\RDset,\WDset,\RRR$ and $\PPP$ of singletons, nonempty finite subsets, nonempty chains, directed subsets, Rudin subsets, well-filtered determined subsets, irreducible subsets and nonempty subsets, respectively.  The familiar inclusions include
\[
 \SSS\leq\CCC\leq\DDD\leq\RDset\leq\WDset\leq\RRR\leq\PPP,
 ~
 \SSS\leq\FFF\leq\PPP;
\]
see~\cite{Xu-Shen-Xi-Zhao-2020}.

In what follows, the capital letter $H$ always denotes a subset system $\HH:\Topz\to\mathbf{Set}$.

\begin{definition}\label{def:H-open}
Let $(X,\tau)$ be a $T_0$-space.  A subset $U\subseteq X$ is $H$-\emph{open relative to} $\tau$ if
\[
  U\cap\cl_\tau A\neq\varnothing~\Longrightarrow~ U\cap A\neq\varnothing
  ~(A\in\HH((X,\tau))).
\]
The family of all such sets is denoted by $\tau_H$.  The topology generated by $\tau_H$ as a subbasis is denoted by $\hdH(\tau)$ and is called the $H$-\emph{derived topology}.  The space $(X,\tau)$ is $H$-\emph{determined} if
\[
     \tau=\hdH(\tau).
\]
The full subcategory of $\Topz$ formed by the $\HH$-determined spaces is denoted by $\DetH$.
\end{definition}

Since every $\tau$-open set satisfies the defining condition, $\tau\subseteq\tau_H$.  Moreover every member of $\tau_H$ is an upper set for the specialization order because $\SSS\leq\HH$.

\begin{proposition}\label{prop:basic-derived}
Let $\tau$ be a $T_0$-topology on $X$.
\begin{enumerate}
\item $\tau\subseteq\tau_H\subseteq\hdH(\tau)\subseteq\alpha(X)$.  In particular, $\tau$ and $\hdH(\tau)$ have the same specialization order.
\item If $\HH_1\leq\HH_2$, then $\tau_{H_2}\subseteq\tau_{H_1}$ and $\hdop{H_2}(\tau)\subseteq\hdop{H_1}(\tau)$.
\item If $\mu\subseteq\nu$ are $T_0$-topologies on the same set, then $\hdH(\mu)\subseteq\hdH(\nu)$.
\item $\hdop{\SSS}(\tau)=\hdop{\FFF}(\tau)=\alpha(X)$ and $\hdop{\PPP}(\tau)=\tau$.
\item The family $\tau_{\CCC}$ is a topology; hence $\hdop{\CCC}(\tau)=\tau_{\CCC}$.
\item If $\HH\leq\DDD$, then $\tau_H$ is a topology and therefore $\hdH(\tau)=\tau_H$.
\end{enumerate}
\end{proposition}

\begin{proof}
For (1), only the last inclusion needs comment.  If $x\in U\in\tau_H$ and $x\leq_X y$, then $x\in U\cap\cl_\tau\{y\}$; since $\{y\}\in\HH(X)$, we obtain $y\in U$.  Thus $U$ is upper.  Since both $\tau$ and $\alpha(X)$ induce the same specialization order, every topology between them does as well.

For (2), increasing the family of test sets strengthens the defining condition.  For (3), the identity $(X,\nu)\to(X,\mu)$ is continuous, so $\HH((X,\nu))\subseteq\HH((X,\mu))$, while $\cl_\nu A\subseteq\cl_\mu A$.

For (4), the closure of a finite set is the finite union of the closures of its points; hence the singleton and finite tests both yield exactly the upper sets.  The assertion for $\PPP$ is the usual closure characterization of open sets.

For (5), arbitrary unions cause no difficulty.  If $U,V\in\tau_C$ and a chain $C$ has closure meeting $U\cap V$, choose $c_1\in C\cap U$ and $c_2\in C\cap V$.  The two points are comparable, and $U,V$ are upper; hence one of $c_1,c_2$ belongs to $U\cap V$.

The same argument proves (6): if $U,V\in\tau_H$ and $A\in\HH(X)\subseteq\DDD(X)$ has closure meeting $U\cap V$, choose $a_1\in A\cap U$ and $a_2\in A\cap V$.  A common upper bound $a\in A$ belongs to $U\cap V$ because both sets are upper.
\end{proof}

\begin{corollary}\label{cor:category-monotone}
If $\HH\leq\HH^*$, every $\HH$-determined space is $\HH^*$-determined.  Thus
\[
   \mathbf{HD}\subseteq\mathbf{H^{*}D}
\]
as full subcategories of $\Topz$.
\end{corollary}

\begin{proof}
If $\tau$ is $\HH$-determined, then
\[
   \tau\subseteq\tau_H\subseteq\hdop{H}(\tau)=\tau,
\]
so $\tau_H=\tau$.  Since $\HH\leq\HH^*$, Proposition~\ref{prop:basic-derived}(2) gives
\[
   \tau\subseteq\tau_{H^*}\subseteq\tau_H=\tau.
\]
Hence $\tau_{H^*}=\tau$, and therefore $\hdop{H^*}(\tau)=\tau$.
\end{proof}

The fixed-point condition also has a useful closed-set formulation.

\begin{proposition}\label{prop:closed-char}
A $T_0$-space $(X,\tau)$ is $\HH$-determined if and only if every subset $C\subseteq X$ satisfying
\[
   A\in\HH(X),\ A\subseteq C
   ~\Longrightarrow~
   \cl_\tau A\subseteq C
\]
is closed.
\end{proposition}

\begin{proof}
Take complements in Definition~\ref{def:H-open}.
\end{proof}

The general construction contains several familiar cases.  For $\HH=\DDD$ one recovers the monotone determined topology of~\cite{Erne-2009, Luo-Xu-2017, Yu-Kou-2015}; for $\HH=\RRR$ one recovers the irreducible derived topology of~\cite{Xu-Yang-Luo-2026}.  At the two extremes, $\SSS$- and $\FFF$-determined spaces are precisely the Alexandroff $T_0$-spaces, while every $T_0$-space is $\PPP$-determined.  Notice already that categorical statements cannot be uniform in $\HH$ without additional hypotheses.

\section{Fixed-point completion and transfinite iteration}

The operator $\hdH$ is monotone and inflationary by Proposition~\ref{prop:basic-derived}, but it need not be idempotent for an arbitrary subset system.  We therefore pass to its least fixed point above a given topology.

\begin{lemma}\label{lem:intersection-fixed}
Let $\{\mu_i:i\in I\}$ be a family of $\HH$-determined topologies on a set $X$, all containing a fixed $T_0$-topology $\tau$.  Then $\mu=\bigcap_{i\in I}\mu_i$ is again an $\HH$-determined $T_0$-topology containing $\tau$.
\end{lemma}

\begin{proof}
The intersection is $T_0$ because it contains the $T_0$-topology $\tau$.  By monotonicity,
\[
   \mu\subseteq\hdH(\mu)\subseteq\hdH(\mu_i)=\mu_i
\]
for every $i$.  Hence $\hdH(\mu)\subseteq\mu$, and equality follows.
\end{proof}

\begin{definition}\label{def:Hd-hull}
For a $T_0$-topology $\tau$ on $X$, define
\[
 \HdH(\tau)=\bigcap\{\mu:\tau\subseteq\mu\text{ and }\mu\text{ is $\HH$-determined}\}.
\]
This family is nonempty because $\alpha(X)$ is $\HH$-determined and $\tau\subseteq \alpha(X)$ (note that $X$ is equipped with the specialization order).
\end{definition}

\begin{theorem}\label{thm:hull-closure}
For every $T_0$-space $(X,\tau)$, we have the following conclusions:
\begin{enumerate}
\item $\HdH(\tau)$ is the least $\HH$-determined topology containing $\tau$;
\item $\tau\subseteq\hdH(\tau)\subseteq\HdH(\tau)\subseteq\alpha(X)$, where $X$ is equipped with the order $\leq_\tau$;
\item $\tau$ and $\HdH(\tau)$ have the same specialization order;
\item the assignment $\tau\mapsto\HdH(\tau)$ is a closure operator on the poset of $T_0$-topologies on $X$: it is inflationary, monotone and idempotent.
\end{enumerate}
\end{theorem}

\begin{proof}
Parts (1) and (2) follow from Lemma~\ref{lem:intersection-fixed}.  Part (3) follows from the two inclusions $\tau\subseteq\HdH(\tau)\subseteq\alpha(X)$.  Monotonicity in (4) follows because the class of $\HH$-determined topologies above a larger initial topology is smaller; idempotence follows from (1).
\end{proof}

We next describe the hull by iteration.  Put $\tau_0=\tau$ and define, for ordinals $\alpha$,
\[
 \tau_{\alpha+1}=\hdH(\tau_\alpha),
 ~
 \tau_\lambda=\bigvee_{\beta<\lambda}\tau_\beta
 ~(\lambda\text{ a limit ordinal}),
\]
where the join is the topology generated by $\bigcup_{\beta<\lambda}\tau_\beta$.

\begin{theorem}\label{thm:transfinite}
There is a least ordinal $\gamma_H(\tau)$ such that
\[
   \tau_{\gamma_H(\tau)+1}=\tau_{\gamma_H(\tau)}.
\]
For this ordinal,
\[
   \HdH(\tau)=\tau_{\gamma_H(\tau)}.
\]
\end{theorem}

\begin{proof}
By transfinite induction,
\[
   \tau\subseteq\tau_\alpha\subseteq\HdH(\tau)\subseteq\alpha(X)
\]
for every $\alpha$, and the chain $(\tau_\alpha)$ is increasing.  Put $\kappa=|\mathcal P(X)|$.  If the iteration did not stabilize before $\kappa^+$, then for every $\alpha<\kappa^+$ we could choose
\[
   U_\alpha\in\tau_{\alpha+1}\setminus\tau_\alpha.
\]
Whenever $\alpha<\beta$, monotonicity gives $U_\alpha\in\tau_\beta$, whereas $U_\beta\notin\tau_\beta$; hence the sets $U_\alpha$ are pairwise distinct.  This would give an injection $\kappa^+\to\mathcal P(X)$, a contradiction.  Thus the iteration stabilizes.

Let $\gamma$ be the first stationary stage.  Then $\tau_\gamma$ is $\HH$-determined and contains $\tau$, so $\HdH(\tau)\subseteq\tau_\gamma$ by minimality.  The reverse inclusion follows from the displayed chain.
\end{proof}

\begin{remark}\label{rem:one-step}
The distinction between $\hdH(\tau)$ and $\HdH(\tau)$ is essential in the general theory.  For several standard systems the first step is already idempotent (for example $\SSS,\FFF,\CCC$ and $\DDD$), while for $\HH=\RRR$ the transfinite construction gives the fixed-point hull considered in~\cite{Xu-Yang-Luo-2026}.  A useful general problem is to characterize those subset systems for which $\hdH$ is idempotent on every $T_0$-topology.
\end{remark}

The construction is functorial.

\begin{proposition}\label{prop:functoriality}
If $f:(X,\tau)\to(Y,\delta)$ is continuous, then the same underlying map is continuous
\[
  f:(X,\hdH(\tau))\longrightarrow(Y,\hdH(\delta))
\]
and
\[
  f:(X,\HdH(\tau))\longrightarrow(Y,\HdH(\delta)).
\]
\end{proposition}

\begin{proof}
Let $V\in\delta_H$ and $A\in\HH((X,\tau))$.  If $f^{-1}(V)\cap\cl_\tau A\neq\varnothing$, then
\[
 V\cap f(\cl_\tau A)\neq\varnothing,
 ~
 f(\cl_\tau A)\subseteq\cl_\delta f(A),
\]
and $f(A)\in\HH((Y,\delta))$.  Hence $V\cap f(A)\neq\varnothing$ and therefore $f^{-1}(V)\in\tau_H$.  This proves the first assertion.  Iterating the argument transfinitely, including the evident limit step, gives continuity at every stage of Theorem~\ref{thm:transfinite} and hence the second assertion.
\end{proof}

For maps whose domain is already a fixed point, continuity has a particularly simple test.

\begin{proposition}\label{prop:continuity-test}
Let $(X,\tau)$ be $\HH$-determined, let $(Y,\delta)$ be a $T_0$-space, and let $f:X\to Y$ be a map.  Then the following are equivalent:
\begin{enumerate}
\item $f:(X,\tau)\to(Y,\delta)$ is continuous;
\item $f(\cl_\tau A)\subseteq\cl_\delta f(A)$ for every $A\in\HH((X,\tau))$;
\item $f:(X,\tau)\to(Y,\HdH(\delta))$ is continuous.
\end{enumerate}
\end{proposition}

\begin{proof}
The implication (1) $\Rightarrow$ (2) is standard, and (1) $\Rightarrow$ (3) follows from Proposition~\ref{prop:functoriality} because $\HdH(\tau)=\tau$.  The implication (3) $\Rightarrow$ (1) follows by composing with the identity $(Y,\HdH(\delta))\to(Y,\delta)$.

For completeness, assume (2).  If $V\in\delta$ and $A\in\HH(X)$ have
$f^{-1}(V)\cap\cl_\tau A\neq\varnothing$, then $V\cap\cl_\delta f(A)\neq\varnothing$.  Since $V$ is open, $V\cap f(A)\neq\varnothing$.  Hence $f^{-1}(V)\in\tau_H=\tau$, proving (1).
\end{proof}

\section{The coreflector and categorical consequences}

We use standard categorical terminology as in, for example, \cite{MacLane-1997}.  Let $I:\DetH\hookrightarrow\Topz$ be the inclusion.  Define
\[
   R_H(X,\tau)=(X,\HdH(\tau)),
   ~
   R_H(f)=f.
\]
Proposition~\ref{prop:functoriality} makes $R_H$ a functor.

\begin{theorem}\label{thm:coreflection}
The functor $R_H:\Topz\to\DetH$ is right adjoint to the inclusion $I$.  Equivalently, $\DetH$ is coreflective in $\Topz$, with coreflection map
\[
   r_X:(X,\HdH(\tau))\longrightarrow(X,\tau),~ r_X(x)=x.
\]
\end{theorem}

\begin{proof}
For $A\in\DetH$ and $X\in\Topz$, Proposition~\ref{prop:continuity-test} gives, naturally in both variables,
\[
   \Topz(IA,X)=\DetH(A,R_HX),
\]
with the same underlying set maps on the two sides.
\end{proof}

This formulation makes the limit theory immediate.

\begin{corollary}\label{cor:complete}
The category $\DetH$ is complete.  If a diagram in $\DetH$ has limit $(L,\lambda)$ in $\Topz$, then its limit in $\DetH$ is
\[
   (L,\HdH(\lambda)).
\]
In particular, for $\HH$-determined spaces $(X_j,\tau_j)$,
\[
   \prod{}^{\DetH}_{j\in J}X_j
   =R_H\!\left(\prod{}^{\Topz}_{j\in J}X_j\right).
\]
\end{corollary}

\begin{proof}
The category $\Topz$ is complete, and the right adjoint $R_H$ preserves limits.
\end{proof}

\begin{proposition}\label{prop:retract}
Every retract of an $\HH$-determined space is $\HH$-determined.
\end{proposition}

\begin{proof}
Let $p:X\to Y$ and $s:Y\to X$ be continuous with $p\circ s=id_Y$, and suppose $X$ is $\HH$-determined.  If $V\in\hdH(\OO(Y))$, Proposition~\ref{prop:functoriality} gives $p^{-1}(V)\in\hdH(\OO(X))=\OO(X)$.  Hence
\[
  V=s^{-1}(p^{-1}(V))\in\OO(Y).
\]
Thus $\hdH(\OO(Y))=\OO(Y)$.
\end{proof}

Reflectivity behaves differently.  We use the criterion of Shen, Xi and Zhao~\cite{Shen-Xi-Zhao-2024}: if a full subcategory $\mathbf K$ of $\Topz$ is closed under homeomorphism and is not contained in $\Topo$, then reflectivity of $\mathbf K$ implies that every sober space belongs to $\mathbf K$.

\begin{proposition}\label{prop:nonreflect-criterion}
For every subset system $\HH$, the category $\DetH$ contains non-$T_1$ spaces.  Consequently, if $\Sob\not\subseteq\DetH$, then $\DetH$ is not reflective in $\Topz$.
\end{proposition}

\begin{proof}
First note that $\DetH$ is closed under homeomorphisms.  Indeed, let $f:(X,\tau)\to(Y,\delta)$ be a homeomorphism and suppose that $(X,\tau)$ is $\HH$-determined.  Proposition~\ref{prop:functoriality} gives a continuous map
\[
   f:(X,\tau)=(X,\hdH(\tau))\longrightarrow(Y,\hdH(\delta)).
\]
Hence $V\in\hdH(\delta)$ implies $f^{-1}(V)\in\tau$, and the original homeomorphism then gives $V\in\delta$.  Thus $\hdH(\delta)\subseteq\delta$, while the reverse inclusion is automatic, so $(Y,\delta)$ is $\HH$-determined.

Every Alexandroff $T_0$-space is $\HH$-determined by Proposition~\ref{prop:basic-derived}(1).  In particular, the Sierpi\'nski space belongs to $\DetH$ and is not $T_1$.  The reflectivity criterion of Shen, Xi and Zhao \cite{Shen-Xi-Zhao-2024} now applies and proves the second statement.
\end{proof}

\begin{theorem}\label{thm:R-nonreflective}
If $\HH$ is an irreducible subset system, i.e. $\HH\leq\RRR$, then $\DetH$ is not reflective in $\Topz$.
\end{theorem}

\begin{proof}
Let $X$ be any non-discrete Hausdorff space.  It is sober, while $\Irr(X)=\SSS(X)$.  Since $\SSS\leq\HH\leq\RRR$, we have $\HH(X)=\SSS(X)$, and therefore
\[
   \hdH(\OO(X))=\alpha(X),
\]
which is the discrete topology because $X$ is $T_1$.  Thus $X$ is not $\HH$-determined, so $\Sob\not\subseteq\DetH$.  Proposition~\ref{prop:nonreflect-criterion} applies.
\end{proof}

\begin{corollary}\label{cor:standard-nonreflective}
The categories of $\SSS$-, $\CCC$-, $\DDD$- and $\RRR$-determined spaces are not reflective in $\Topz$.  The same holds for $\FFF$-determined spaces because the $\SSS$- and $\FFF$-determined spaces coincide.
\end{corollary}

\begin{remark}\label{rem:P-reflective}
The hypothesis in Theorem~\ref{thm:R-nonreflective} cannot simply be omitted.  For $\HH=\PPP$, Proposition~\ref{prop:basic-derived}(4) gives $\hdH(\tau)=\tau$ for every $T_0$-topology, so $\mathbf{HD}=\Topz$; in particular this category is reflective as well as coreflective.  This also shows why a blanket non-reflectivity assertion for arbitrary subset systems is false.
\end{remark}

\section{Subspaces: stability hypotheses on the subset system}

For the irreducible system $\RRR$, irreducibility of a subset is unchanged when it is regarded inside a subspace containing it; this fact underlies the hereditary results in~\cite{Xu-Yang-Luo-2026}.  An arbitrary subset system need not have this property, and functoriality of $\HH$ only gives one direction.  We therefore make the required hypotheses explicit.

\begin{definition}\label{def:restriction-local}
A subset system $\HH$ is called
\begin{enumerate}
\item[(1)] \emph{restriction-stable} if, whenever $Y$ is a nonempty subspace of $X$ and $A\subseteq Y$, the implication
\[
   A\in\HH(X)~\Longrightarrow~ A\in\HH(Y)
\]
holds;
\item[(2)] \emph{open-local} if, whenever $U$ is open in $X$ and $A\in\HH(X)$ with $A\cap U\neq\varnothing$, one has
\[
   A\cap U\in\HH(U).
\]
\end{enumerate}
\end{definition}

The reverse implication in (1) is automatic from the subset-system axioms because the inclusion $Y\hookrightarrow X$ is continuous.  The standard systems $\SSS,\FFF,\CCC,\DDD$ and $\RRR$ are restriction-stable and open-local.  For $\DDD$, for example, a nonempty intersection of a directed set with an open upper set is again directed; for $\RRR$, a nonempty intersection of an irreducible set with an open set is irreducible in the open subspace.

\begin{theorem}\label{thm:closed-hereditary}
If $\HH$ is restriction-stable, then every nonempty closed subspace of an $\HH$-determined space is $\HH$-determined.
\end{theorem}

\begin{proof}
Let $X$ be $\HH$-determined, let $C\subseteq X$ be nonempty and closed, and let $V\subseteq C$ be $\HH$-open relative to the subspace topology.  Put $U=V\cup(X\setminus C)$.  We claim that $U\in\OO(X)$.

Suppose $A\in\HH(X)$ and $U\cap\cl_XA\neq\varnothing$.  If $A\cap U\neq\varnothing$ there is nothing to prove.  Otherwise $A\subseteq C\setminus V$.  Restriction-stability gives $A\in\HH(C)$, while $C$ is closed, so $\cl_XA=\cl_CA\subseteq C$.  Hence
\[
   \varnothing\neq U\cap\cl_XA=V\cap\cl_CA,
\]
and the defining property of $V$ yields $A\cap V\neq\varnothing$, a contradiction.  Thus $U$ is $\HH$-open in $X$ and therefore open because $X$ is $\HH$-determined.  Consequently $V=U\cap C$ is open in $C$.
\end{proof}

\begin{theorem}\label{thm:open-hereditary}
If $\HH$ is open-local, then every nonempty open subspace of an $\HH$-determined space is $\HH$-determined.
\end{theorem}

\begin{proof}
Let $X$ be $\HH$-determined, $O\subseteq X$ be nonempty and open, and let $V\subseteq O$ be $\HH$-open relative to $O$.  We show directly that $V$ is $\HH$-open in $X$.

Take $A\in\HH(X)$ with $V\cap\cl_XA\neq\varnothing$.  Then $O\cap\cl_XA\neq\varnothing$, and since $O$ is open, $A\cap O\neq\varnothing$.  Put $B=A\cap O$.  By open-locality, $B\in\HH(O)$.  Moreover
\[
 A\subseteq B\cup(X\setminus O)
\]
and $X\setminus O$ is closed, so
\[
  \cl_XA\subseteq\cl_XB\cup(X\setminus O).
\]
Intersecting with $V\subseteq O$ gives $V\cap\cl_XB\neq\varnothing$, equivalently $V\cap\cl_OB\neq\varnothing$.  Hence $V\cap B\neq\varnothing$, and therefore $V\cap A\neq\varnothing$.  Thus $V$ is $\HH$-open in $X$, so it is open.  In particular it is open in $O$.
\end{proof}

\begin{corollary}\label{cor:standard-hereditary}
For $\HH\in\{\SSS,\FFF,\CCC,\DDD,\RRR\}$, the property of being $\HH$-determined is both closed-hereditary and open-hereditary.
\end{corollary}

For $\HH=\RRR$, Corollary~\ref{cor:standard-hereditary} recovers the corresponding results of~\cite{Xu-Yang-Luo-2026}; the point here is that the general argument identifies precisely the properties of the subset system that are used.

\section{Cartesian closedness below the directed subset system}

Assume throughout this section that
\[
   \HH\leq\DDD.
\]
For $X,Y\in\DetH$, write $X\times_HY$ for their categorical product in $\DetH$, namely the $\HH$-determined hull of their product in $\Topz$ as in Corollary~\ref{cor:complete}.

For $T_0$-spaces $Y,Z$, let $C(Y,Z)$ be the set of continuous maps from $Y$ to $Z$.  Equip $C(Y,Z)$ with the topology $\ptop(Y,Z)$ of pointwise convergence, whose subbasic open sets are
\[
   [y,U]=\{f\in C(Y,Z):f(y)\in U\},
   ~ y\in Y,\ U\in\OO(Z).
\]
The pointwise topology is $T_0$.  Indeed, if $f\neq g$, choose $y\in Y$ with $f(y)\neq g(y)$ and use the $T_0$ property of $Z$ to find an open set that contains exactly one of $f(y)$ and $g(y)$; the corresponding subbasic set separates $f$ and $g$.  Moreover, its specialization order is the pointwise order: $f\leq g$ if and only if $f(y)\leq_Z g(y)$ for every $y\in Y$.  Define
\[
   [Y,Z]_H:=R_H(C(Y,Z),\ptop(Y,Z)).
\]
Thus $[Y,Z]_H$ is automatically an object of $\DetH$.

\begin{lemma}\label{lem:evaluation}
For $Y,Z\in\DetH$, the evaluation map
\[
  \operatorname{ev}:[Y,Z]_H\times_HY\longrightarrow Z,
  ~ \operatorname{ev}(f,y)=f(y),
\]
is continuous.
\end{lemma}

\begin{proof}
We use Proposition~\ref{prop:continuity-test}.  Let
\[
   A\in\HH([Y,Z]_H\times_HY)
\]
and let $(f,y)\in\cl_{[Y,Z]_H\times_HY}A$.  The two projections are continuous, hence
\[
   f\in\cl_{[Y,Z]_H}p_1(A),
   ~
   y\in\cl_Yp_2(A).
\]
The topology of $[Y,Z]_H$ is finer than the pointwise topology, so $f$ also belongs to the pointwise closure of $p_1(A)$.

Let $U\in\OO(Z)$ contain $f(y)$.  Since $f:Y\to Z$ is continuous, $f^{-1}(U)$ is an open neighborhood of $y$.  Choose $(g_1,y_1)\in A$ with $y_1\in f^{-1}(U)$.  Then $f(y_1)\in U$.  Because $f$ lies in the pointwise closure of $p_1(A)$, the subbasic neighborhood $[y_1,U]$ meets $p_1(A)$; choose $(g_2,y_2)\in A$ such that $g_2(y_1)\in U$.

Now $A\in\HH([Y,Z]_H\times_HY)\subseteq\DDD([Y,Z]_H\times_HY)$, so $A$ is directed in the specialization order.  Choose $(g_0,y_0)\in A$ above both $(g_1,y_1)$ and $(g_2,y_2)$.  By Theorem~\ref{thm:hull-closure}(3), the categorical product has the same specialization order as the ordinary product, and the specialization order on $[Y,Z]_H$ is the pointwise order.  Since every continuous map is order preserving,
\[
   g_2(y_1)\leq g_2(y_0)\leq g_0(y_0).
\]
Because $U$ is an upper set, $g_0(y_0)\in U$.  Thus every open neighborhood of $f(y)$ meets $\operatorname{ev}(A)$, and hence
\[
   \operatorname{ev}\bigl(\cl_{[Y,Z]_H\times_HY}A\bigr)
   \subseteq
   \cl_Z\operatorname{ev}(A).
\]
Proposition~\ref{prop:continuity-test} now yields continuity of $\operatorname{ev}$.
\end{proof}

\begin{theorem}\label{thm:ccc}
If $\HH\leq\DDD$, then $\DetH$ is Cartesian closed.  For $Y,Z\in\DetH$, the object $[Y,Z]_H$ defined above is an exponential of $Z$ by $Y$.
\end{theorem}

\begin{proof}
Let $X,Y,Z\in\DetH$.  Given a continuous map
\[
    q:X\times_HY\longrightarrow Z,
\]
define $\widehat q:X\to C(Y,Z)$ by $\widehat q(x)(y)=q(x,y)$.  For each fixed $x$, the section $y\mapsto(x,y)$ is the pairing of the constant map $Y\to X$ with $\operatorname{id}_Y$ and is therefore continuous by the universal property of the categorical product.  Thus $\widehat q(x):Y\to Z$ is continuous.  Similarly, for each fixed $y$, the section $x\mapsto(x,y)$ is continuous, so $x\mapsto q(x,y)$ is continuous.  Hence $\widehat q:X\to(C(Y,Z),\ptop(Y,Z))$ is continuous by the defining initial property of the pointwise topology.  Since $X$ is $\HH$-determined, Theorem~\ref{thm:coreflection} (or Proposition~\ref{prop:continuity-test}) lifts this uniquely to a continuous map
\[
   \widehat q:X\longrightarrow[Y,Z]_H.
\]

Conversely, if $g:X\to[Y,Z]_H$ is continuous, then
\[
  X\times_HY\xrightarrow{\ g\times\operatorname{id}_Y\ }[Y,Z]_H\times_HY
  \xrightarrow{\ \operatorname{ev}\ } Z
\]
is continuous by Lemma~\ref{lem:evaluation}.  These two constructions are inverse to one another on underlying maps.  Therefore there is a natural bijection
\[
   \DetH(X\times_HY,Z)\cong\DetH(X,[Y,Z]_H),
\]
which is precisely the exponential universal property.
\end{proof}

\begin{corollary}\label{cor:ccc-standard}
The categories of $\SSS$-, $\FFF$- and $\CCC$-determined spaces are Cartesian closed.
\end{corollary}

\begin{proof}
The $\SSS$- and $\FFF$-determined spaces coincide by
Proposition~\ref{prop:basic-derived}(4), while $\SSS,\CCC\leq\DDD$.
Apply Theorem~\ref{thm:ccc}.
\end{proof}

\begin{corollary}\label{cor:ddd-standard}
The category of monotone determined spaces is Cartesian closed; see also \cite{Yu-Kou-2015,Luo-Xu-2017}.
\end{corollary}

\section{Concluding remarks and questions}

The fixed-point formulation separates three levels of the construction: the one-step topology $\hdH(\tau)$, its idempotent hull $\HdH(\tau)$, and the categorical coreflector $R_H$. This makes clear which conclusions are formal consequences of the subset-system axioms and which require additional structure on $\HH$.

From the standpoint of programming semantics, the main categorical consequence is Theorem~\ref{thm:ccc}: for $\HH\leq\DDD$, the same fixed-point mechanism that produces the coreflective category also supplies its exponentials. The paper does not commit to a particular programming language or type system; rather, it isolates a reusable topological and categorical mechanism that can be tested in future against concrete domain-theoretic models and higher-order type interpretations.

The following questions arise naturally.

\begin{question}
Characterize the subset systems $\HH$ for which $\hdH$ is idempotent on every $T_0$-topology, equivalently, for which $\hdH(\tau)=\HdH(\tau)$ always holds.
\end{question}

\begin{question}
Can the restriction-stability and open-locality assumptions in Section~5 be replaced by weaker intrinsic conditions that are still necessary and sufficient for closed- and open-heredity of $\HH$-determined spaces?
\end{question}

\begin{question}
Characterize the subset systems $\HH$ for which $\DetH$ is reflective in $\Topz$.  Theorem~\ref{thm:R-nonreflective} gives a broad negative class, while $\HH=\PPP$ gives the opposite extreme.
\end{question}

\begin{question}
Can the hypothesis $\HH\leq\DDD$ in Theorem~\ref{thm:ccc} be weakened?  More generally, which properties of $\HH$ guarantee continuity of the evaluation map for the coreflected pointwise function space?
\end{question}

\backmatter

\section*{Statements and Declarations}

\bmhead{Funding}
This research was supported by the National Natural Science Foundation of
China (Nos.~12471070 and 12071199).

\bmhead{Author contributions}
Xiaoquan Xu: Conceptualization,  Investigation, methodology,  formal analysis, writing-review and editing.
Xinjie Huang:  Conceptualization, investigation,  formal analysis, writing-original draft.
Yi Yang:  Investigation,  formal analysis, writing-review and editing
All authors have reviewed and approved the final manuscript.


\bmhead{Competing interests}
The authors have no relevant financial or non-financial interests to disclose.

\bmhead{Data availability}
No datasets were generated or analysed in this theoretical study.

\bmhead{Use of AI-assisted tools}
During the preparation of this work, the authors used ChatGPT (OpenAI) to assist with organization, language and readability.
The authors reviewed and edited the resulting text, verified the mathematical
content and references, and take full responsibility for the final manuscript.

\end{document}